\documentclass[a4paper,12pt, twoside]{article}
\usepackage{ThPackage1}
\usepackage{ThPackage2}

\begin{document}

%
%
%
%
%
%

\begin{center}\textbf{\Large The unit theorem for finite-dimensional algebras} \end{center}
\begin{center} Hendrik Lenstra \end{center}
Notes for a mini-course given at the Max Planck Institute for Mathematics in Bonn, on October 25th and November 3rd 2016.\\
\\
\textit{Summary.} The \textit{unit theorem} to which the present mini-course is devoted is a theorem from algebra that has a combinatorial flavour, and that originated in fact from algebraic combinatorics. Beyond a proof, the course also addresses applications, one of which is a proof of the normal basis theorem from Galois theory.
\vspace{\baselineskip}

\tableofcontents

\subsection{Introduction}
An \textit{algebra} over a commutative ring $F$ is a ring $A$ equipped with a ring homomorphism $F \lra \Ze(A) = (\textnormal{center of } A)$ (so that the multiplication map $A\times A \lra A$ is $F$-bilinear). Let $F$ be a field. An algebra $A$ over $F$ is called \textit{finite-dimensional} if the dimension of $A$, viewed as a vector space over $F$, is finite.

\begin{thm}[Unit theorem]
Let $F$ be a field, let $A$ be a finite-dimensional algebra over $F$, and let $H$ be an additive subgroup of $A$ that spans $A$ as an $F$-vector space. Then $H$ contains a unit of $A$.
\end{thm}

\noindent If $H$ is as in the theorem, then $H$ contains some \textit{basis} $B$ for $A$ as a vector space over $F$, and $H$ will then also contain the additive subgroup of $A$ generated by $B$; that additive subgroup equals $\sum_{b\in B} (\Z\cdot 1)\cdot b$, where $\Z\cdot 1$ denotes the prime subring of $F$ (so $\Z\cdot 1 \cong \Z/(\Char F) \Z)$. Conversely, if $B$ is any $F$-vector space basis for $A$, then $\sum_{b\in B}(\Z\cdot 1)\cdot b$ is an additive subgroup of $A$ that spans $A$ as an $F$-vector space. So we see that the unit theorem is equivalent to the following statement: if $A$ is a finite-dimensional algebra over a field $F$, with vector space basis $B$, then $A$ contains a unit all of whose coefficients on $B$ belong to the prime subring of $F$.

\begin{exe}
Show that the theorem is incorrect without the finite dimensionality assumption.
\end{exe}

\noindent Interesting examples of finite-dimensional algebras over $F$ are: matrix rings over $F$; group rings of finite groups over $F$; extension fields of finite degree over $F$ (for which the unit theorem is trivial); and rings of the form $F^n = F \times F \times \ldots \times F$ ($n$ copies of $F$, for some $n \in \Z_{\ge 0})$, with componentwise ring operations. Algebras of the latter type are said to be \textit{totally split}. For a totally split algebra, the unit theorem says: if $B$ is a basis for $F^n$ over $F$, then some \textit{integral} linear combination of $B$ has all coordinates different from zero.
\par
The three main subjects of the present mini-course are a \textit{proof} of the unit theorem, \textit{applications}, and the case of \textit{commutative} algebras.

\begin{theproof} To explain the main idea and the main difficulty of the proof, I give it in a special case. Let $F$ and $A$ be as in the unit theorem. As we saw above, we may assume that $H$ is the additive subgroup of $A$ generated by some vector space basis $b_1,b_2,\ldots,b_n$ of $A$ over $F$. For $x = \sum_{i=1}^n x_ib_i \in A$ (with $x_i \in F$), let $M(x)$ be the $n\times n$-matrix that describes the $F$-linear map $A \lra A$, $a \mapsto xa$, on the basis $b_1,b_2,\ldots,b_n$. Then $x$ is a unit of $A$ if and only if $\det M(x)\neq 0$. Also, one has $M(x) = \sum_{i=1}^n x_iM(b_i)$, so if one defines the element $P$ of the polynomial ring $F[X_1,\ldots,X_n]$ by \[P = \det\mleft(\sum_{i=1}^n X_iM(b_i)\mright)\] then one has $\det M(x) = P(x_1,x_2,\ldots,x_n)$, and $x$ is a unit of $A$ if and only if $P(x_1,x_2,\ldots,x_n) \neq 0$. Since for $x = 1$ one has $\det M(x) = 1 \neq 0$, the polynomial $P$ is not the zero polynomial. Now make the assumption that $F$ has \textit{characteristic zero} or, equivalently, that its prime subring is infinite. We may identify that prime subring with $\Z$. Since, by a well-known theorem (see \S 2), a non-zero polynomial in $n$ variables over a field does not vanish identically on any product of $n$ infinite subsets of the field, there exists $(x_1,x_2,\ldots,x_n) \in \Z\times \Z \times \ldots \times \Z$ with $P(x_1,x_2,\ldots,x_n) \neq 0$. This gives rise to an element $x = \sum_{i=1}^n x_ib_i$ that belongs to $H$ and is a unit, as desired.

This proves the unit theorem in the case of characteristic \textit{zero}. In the case of non-zero characteristic $p$, the main difficulty is that a non-zero polynomial may vanish identically on a product of finite sets, such as $\Z/p\Z \times \Z/p\Z \times \ldots \times \Z/p\Z$. To overcome this difficulty, we shall develop a variant of the \textit{combinatorial Nullstellensatz} that is independent of the choice of coordinates; it states that, under certain conditions, the phenomenon just described does not occur. To prove that in our case those conditions are satisfied, we use a classical structure theorem for finite-dimensional algebras over an algebraically closed field.
\end{theproof}
\begin{exe}
Let $F$, $A$, $H$ be as in the unit theorem, and assume that $F$ has characteristic zero. Prove that every coset $a+H$ of $H$, with $a \in A$, contains a unit of $A$. 
\end{exe}

\begin{applications} 
To illustrate the applications of the unit theorem, we use it to obtain a mild sharpening of the \textit{normal basis theorem} from Galois theory. Let $E \subset F$ be a finite Galois extension of fields, with Galois group $G$. It is well-known that one has $\# G = [F:E]$. The normal basis theorem asserts that there exists $\alpha \in F$ such that $(\sigma \alpha)_{\sigma \in G}$ is a basis for $F$ as a vector space over $E$. We shall prove the following.
\end{applications}

\begin{thm}
Let $E$, $F$, $G$ be as above, and let $H$ be an additive subgroup of $F$ that spans $F$ as an $E$-vector space. Then there exists $\alpha \in H$ such that $(\sigma \alpha)_{\sigma \in G}$ is a vector space basis of $F$ over $E$.
\end{thm}

\begin{exe}
Verify this theorem directly in the case $E = \R$, $F=\C$.
\end{exe}

\noindent To deduce the theorem from the unit theorem, we let $F[G]$ denote the group ring of $G$ over $F$, and we define $\varphi\colon F \lra F[G]$ by
\[ \varphi(\alpha) = \sum_{\tau \in G} (\tau^{-1} \alpha)\cdot \tau.\]
This map has the following properties:
\begin{itemize}
\item[(i)] $\varphi$ is $E$-linear,
\item[(ii)] $\varphi(\sigma \alpha) = \sigma \cdot \varphi(\alpha)$ for all $\sigma \in G$, $\alpha \in F$,
\item[(iii)] for any basis $(b_i)_{i\in I}$ for $F$ as an $E$-vector space, $(\varphi(b_i))_{i\in I}$ is a basis for $F[G]$ as an $F$-vector space.
\end{itemize}
Here (i) is obvious, and (ii) follows from an easy computation. For (iii), it suffices to show that the square matrix $(\tau b_i)_{\tau \in G, i \in I}$ is invertible. Suppose $a_\tau \in F$, for $\tau \in G$, are such that for all $i\in I$ one has $\sum_{\tau \in G} a_\tau \tau b_i = 0$. Then for all $b\in F$ one has $\sum_{\tau \in G} a_\tau \tau b = 0$, so by Artin's theorem on linear independence of field homomorphisms one has $a_\tau = 0$ for all $\tau \in G$, and (iii) follows.

Now let $H$ be as in the theorem. Then $H$ spans $F$ as an $E$-vector space, so (iii) shows that $\varphi(H)$ spans $F[G]$ as an $F$-vector space. Applying the unit theorem, with $F[G]$ and $\varphi(H)$ in the roles of $A$ and $H$, we can choose $\alpha \in H$ such that $\varphi(\alpha)$ is a unit of $F[G]$. Then $F[G]\cdot \varphi(\alpha) = F[G]$, so $(\sigma \cdot \varphi(\alpha))_{\sigma \in G}$ is linearly independent over $F$. Hence (ii) shows that $(\varphi(\sigma \alpha))_{\sigma \in G}$ is linearly independent over $F$, so \textit{a fortiori} over $E$. From (i) it now follows that $(\sigma \alpha)_{\sigma \in G}$ is linearly independent over $E$, so it is a basis for $F$ over $E$, as required.

\begin{exe}
Let $E$, $F$, $G$, $\varphi$ be as above.
\begin{itemize}[noitemsep,nolistsep]
\item[(a)] Prove that $\varphi$ is a ring homomorphism if and only if $E = F$.
\item[(b)] Let $\alpha \in F$. Prove: $(\sigma \alpha)_{\sigma \in G}$ is a basis for $F$ over $E$ if \textit{and only if} $\varphi(\alpha)$ is a unit of $F[G]$.
\end{itemize}
\end{exe}

\noindent The normal basis theorem is equivalent to the assertion that $F$, when viewed as a module over the group ring $E[G]$, is isomorphic to $E[G]$ itself. Paraphrasing the proof just given, one may say that it starts by verifying this assertion ``after base extension from $E$ to $F$''; more precisely, the properties of $\varphi$ show that there is an isomorphism $F\otimes_E F \cong E[G] \otimes_E F$ of modules over the ring $E[G]\otimes_E F$, which is nothing but the group ring $F[G]$. Next, the unit theorem is used to ``descend'' this isomorphism from $F$ to $E$. All applications of the unit theorem that we shall give, are to similar \textit{descent problems}.

\begin{commcase}In the commutative case, we start with the following exercise, in which $R^\ast$ denotes the group of units of a ring $R$.\end{commcase}

\begin{exe}
Let $A$ be a commutative finite-dimensional algebra over some field $F$. Prove that for some $n \in \Z_{\ge 0}$ there exists a vector space isomorphism $\psi\colon F^n \lra A$ such that $\psi(F^{\ast n}) \subset A^\ast$.
\end{exe}

\noindent This exercise shows that, in order to prove the unit theorem in the case $A$ is commutative, one may assume that $A$ is \textit{totally split}, in which case, as we saw above, the unit theorem allows a reformulation that does not refer to the ring structure.

In the commutative case, the unit theorem admits several different proofs, each of which supplies additional information. For example, we have the following result, which I learnt from Bas Edixhoven.

\begin{thm}
Let $F$, $A$, $H$ be as in the unit theorem, and write $n$ for the dimension of $A$ as an $F$-vector space. Assume that the characteristic $p$ of $F$ is non-zero, and that $A$ is commutative. Then the number of units of $A$ contained in $H$ is at least $(p-1)^n$.
\end{thm}

\noindent One may wonder whether the theorem remains true without the commutativity assumption.

\subsection{The combinatorial Nullstellensatz}

\noindent Throughout this section, $F$ denotes a field. For a non-zero polynomial $g$ in several variables, we write $\deg g$ for the total degree of $g$.

In the introduction we used the well-known fact that if $g\in F[X_1,\ldots,X_n]$ is non-zero, with $n\in\Z_{\ge 0}$, and $S_1,S_2,\ldots, S_n$ are infinite subsets of $F$, then $g$ does not vanish identically on $S_1\times S_2 \times \ldots \times S_n$. The \textit{combinatorial Nullstellensatz} is a quantitative refinement of this statement.

\begin{thm}[Combinatorial Nullstellensatz, Noga Alon, 1999] Let $F$ be a field, let $n\in \Z_{\ge 0}$, $d_1$, $d_2$, $\ldots$, $d_n \in \Z_{\ge 0}$, and let $g\in F[X_1,\ldots,X_n]$. Suppose that $g$ has a non-zero coefficient at $X_1^{d_1}\cdots X_n^{d_n}$ and that $d_1+d_2+\ldots+d_n = \deg g$. Let $S_i \subset F$ satisfy $\# S_i > d_i$, for $1 \le i \le n$. Then $g$ does not vanish identically on $S_1 \times \ldots \times S_n$.
\end{thm}

\begin{prfd}
For $x = (x_1,x_2,\ldots,x_n) \in F^n$ the map $F[X_1,\ldots,X_n] \lra F$, $f\mapsto f(x) = f(x_1,x_2,\ldots,x_n)$ is a surjective ring homomorphism, so denoting its kernel by $I(x)$, we have $F[X_1,\ldots,X_n]/I(x) \cong F$. We have 
\[ I(x) + I(y) = F[X_1,\ldots,X_n] \textnormal{ for } x = (x_j)_{j=1}^n, y= (y_j)_{j=1}^n \in F^n, x \neq y,\]
since if $x_j \neq y_j$ then $I(x) + I(y)$ contains $(X_j-x_j) - (X_j - y_j) = y_j - x_j$, which is a unit. The Chinese remainder theorem now implies that, for any finite subset $S \subset F^n$, the map 
\[F[X_1,\ldots,X_n] \lra \prod_{x\in S} F[X_1,\ldots,X_n]/I(x) \cong F^S\]
sending $f \in F[X_1,\ldots, X_n]$ to $(f(x))_{x \in S}$ in $F^S$, is a \textit{surjective} ring homomorphism.

Turning to the proof of the theorem, we assume, without loss of generality, that each $S_i$ is \textit{finite} of cardinality $d_i + 1$, and we apply the above to $S = S_1 \times S_2 \times \ldots \times S_n$. For each $i$, the polynomial
\[ h_i = \prod_{x\in S_i} (X_i - x)\]
is a monic polynomial in $X_i$ of degree $d_i + 1$, and it belongs to the kernel of the ring homomorphism $F[X_1,X_2,\ldots, X_n] \lra F^S$ that we defined. We claim that this ring homomorphism is still surjective when restricted to the $F$-vector space $V$ (say) spanned by
\[X_1^{i_1}X_2^{i_2} \cdots X_n^{i_n} \textnormal{ with all } i_j \in \{0,1,\ldots, d_j\}.\]
To see this, it suffices to observe that for each $f \in F[X_1,X_2,\ldots,X_n]$ one can perform successive ``divisions with remainder'' by $h_1,h_2,\ldots,h_n$ to find $\ri(f) \in V$ with $f\equiv \ri(f) \textnormal{ mod }(h_1,h_2,\ldots,h_n)$, so that $f$ and $\ri(f)$ have the same image in $F^S$.

The hypothesis $d_1+d_2+\ldots+d_n = \deg g$ ensures that $g$ and $\ri(g)$ have the same coefficient at $X_1^{d_1}\cdots X_n^{d_n}$. Since it was supposed to be non-zero, it follows that $\ri(g) \neq 0$.

Because $V$ and $F^S$ both have $F$-dimension $\prod_{j=1}^n (d_j+1)$, the surjective map $V\lra F^S$, which is $F$-linear, is actually an isomorphism of $F$-vector spaces. Hence $\ri(g)$ does not map to $0$. It follows that $g$ does not map to $0 \in F^S$, which is the statement to be proved. This proves the combinatorial Nullstellensatz.
\end{prfd}

\begin{exe}
Prove that, in the situation of the proof just given, the kernel of the surjective ring homomorphism $F[X_1,\ldots,X_n] \lra F^S$ is for $S = S_1 \times S_2 \times \ldots \times S_n$ generated by $h_1,h_2,\ldots, h_n$. $^\ast$Can you also prove that for \textit{any} finite subset $S \subset F^n$ the kernel can be generated by $n$ elements?
\end{exe}

\begin{exe}
Prove that the combinatorial Nullstellensatz remains valid if the conditions on $g$ are replaced by the following: $g$ has a non-zero term $cX_1^{d_1} \cdots X_n^{d_n}$ (with $c\in F$) such that no other non-zero term $c' X_1^{e_1} \cdots X_n^{e_n}$ (with $c' \in F$) of $g$ satisfies $e_1\ge d_1,\ldots, e_n \ge d_n$.
\end{exe}

\noindent The combinatorial Nullstellensatz has seen many dramatic applications in extremal graph theory and in arithmetic combinatorics. As an example of the latter type, we give a quick proof of the \textit{theorem of Cauchy--Davenport}. It has seen many proofs through the ages, but none is as brief as the present one.

\begin{thm}[Cauchy, 1813; Davenport, 1935]
Let $p$ be a prime number, let $A$, $B \subset \Z/p\Z$ be non-empty subsets, and put $A+B = \{ a+b : a \in A, b \in B\}$. Then we have $\# (A+B) \ge \min\{\# A + \# B - 1, p\}$.
\end{thm}
\begin{prfd}
If $\# A + \# B > p$, then for each $c \in \Z/p\Z$ the sets $A$ and $c - B$ must intersect, and therefore $A+B = \Z/p\Z$. Assume therefore that $\#A + \# B \le p$. If the claim fails, then we can choose $C \subset\Z/p\Z$ with $A+B \subset C$ and $\# C = \# A + \# B - 2$. Then $g = \prod_{c \in C} (X_1 + X_2 - c)$ vanishes on $A \times B$, has degree $\# A + \# B - 2$, and has a non-zero coefficient at $X_1^{\# A - 1} \cdot X_2^{\#B - 1}$, contradicting the combinatorial Nullstellensatz. This proves the theorem.
\end{prfd}

\noindent We next prove that, in a special case, the set of $g\in F[X_1,\ldots,X_n]$ to which the combinatorial Nullstellensatz applies, is invariant under invertible linear substitutions of $X_1,\ldots, X_n$.

\begin{lem}
Suppose that the characteristic $p$ of the field $F$ is non-zero, and let $m$, $n \in \Z_{\ge 0}$. Denote by $G$ the group of ring automorphisms of $F[X_1,\ldots,X_n]$ that are the identity on $F$ and that map $F\cdot X_1 + \ldots + F\cdot X_n$ to itself, and write $\D$ for the set of all $g \in F[X_1,\ldots,X_n]$ that have a non-zero term $cX_1^{d_1}\cdots X_n^{d_n}$ (with $c\in F$) satisfying
\[ d_i < p^m \textit{ for all $i$, and } \sum_{i=1}^n d_i = \deg g\textit{.}\]
Then for all $\sigma \in G$ we have $\sigma\D = \D$.
\end{lem}

\noindent Note that $\D$ consists of those polynomials that by the combinatorial Nullstellensatz are guaranteed \textit{not} to vanish identically on any set of the form $S_1 \times S_2 \times \ldots \times S_n$, where each $S_i \subset F$ satisfies $\# S_i = p^m$.

\begin{prfd}[Proof of the lemma]
For $g\in F[X_1,\ldots,X_n]$, we write $\lt g$ for the sum of the terms of degree $\deg g$ of $g$ if $g\neq 0$, and $\lt 0 = 0$. We have $g\notin \D$ if and only if every term of $\lt g$ is for some $j$ divisible by $X_j^{p^m}$. In other words, writing $I$ for the ideal $(X_1^{p^m},\ldots, X_n^{p^m})$ of $F[X_1,\ldots,X_n]$, we have
\[ g\notin \D \iff \lt g \in I.\]
Let $\sigma \in G$, and put $V = F\cdot X_1 + \ldots + F\cdot X_n$. Since the map $v\mapsto v^{p^m}$ is additive, $I$ is the same as the ideal generated by $\{v^{p^m} : v \in V\}$. Thus, from $\sigma V = V$ it follows that $\sigma I = I$. It also follows that $\sigma$ preserves degree and homogeneity of polynomials, and therefore that $\lt \sigma(g) = \sigma(\lt g)$. Thus we have
\[ g \notin \D \iff \lt g \in I \iff \sigma(\lt g) \in \sigma I \iff \lt(\sigma g) \in I \iff \sigma g \notin \D.\]
This implies the lemma.
\end{prfd}

\begin{exe}
Let $F$, $p$, $m$, $n$, $G$, $\D$ be as in the lemma. 
\begin{itemize}[noitemsep,nolistsep]
\item[\rm(a)] Prove: $G$ is isomorphic to the group $\GL(n,F)$ of invertible $n\times n$-matrices over $F$.
\item[\rm(b)] Let $\Gamma$ be the group of all ring automorphisms $\sigma$ of $F[X_1,\ldots,X_n]$ with the property that $\deg(\sigma g) = \deg g$ for all non-zero $g\in F[X_1,\ldots,X_n]$. Prove: for all $\sigma \in \Gamma$ one has $\sigma \D = \D$.
\item[\rm(c)] Prove that the group $\Gamma$ from (b) is isomorphic to a repeated semidirect product $F^n \rtimes (\GL(n,F) \rtimes \Aut F)$.
\end{itemize}
\end{exe}

\begin{thm}[$p$-power non-vanishing theorem] 
Let $F$ be a field of non-zero characteristic $p$, let $m$, $n \in \Z_{\ge 0}$, $d_1$, $d_2$, $\ldots$, $d_n \in \Z_{\ge 0}$, and suppose $d_i < p^m$ for all $i$. Let $g \in F[X_1,\ldots,X_n]$ be a polynomial with $\deg g = d_1 + d_2 + \ldots + d_n$ that has a non-zero coefficient at $X_1^{d_1}\cdots X_n^{d_n}$. Let $S_i \subset F$ satisfy $\# S_i = p^m$, for $1 \leq i \leq n$, and let $b_1,b_2,\ldots, b_n$ be a basis for the $F$-vector space $F^n$. Then $g$ does not vanish identically on $\sum_{i=1}^n S_ib_i = \{\sum_{i=1}^n x_i b_i : x_i \in S_i\ (1\le i \le n)\}$.
\end{thm}
\begin{prfd}
Write $b_i = (b_{ij})_{j=1}^n$. Then the matrix $(b_{ij})_{i,j=1}^n$ belongs to $\GL(n,F)$, and the map \[f\mapsto \sigma f = f\mleft(\sum_i X_ib_{i1},\ldots,\sum_i X_ib_{in}\mright)\] is a ring automorphism $\sigma$ of $F[X_1,\ldots,X_n]$ that belongs to the group $G$ defined in the lemma. With $\D$ as in the lemma, we have $g \in \D$, so by the lemma we have $\sigma g \in \D$. By the combinatorial Nullstellensatz, there are $x_i \in S_i$ (for $1\le i \le n$) with $(\sigma g)(x_1,\ldots,x_n) \neq 0$. But \[(\sigma g)(x_1,\ldots,x_n) = g\mleft(\sum_i x_ib_{i1},\ldots, \sum_i x_i b_{in}\mright) = g\mleft(\sum_i x_i b_i\mright),\] and the theorem follows.
\end{prfd}

\begin{thm}[Non-vanishing on subgroups] 
Let $F$ be a field of non-zero characteristic $p$, let $n \in \Z_{\ge 0}$, and let $g \in F[X_1,\ldots, X_n]$. Suppose that there exist integers $d_1$, $\ldots$, $d_n \in \{0,1,\ldots, p-1\}$ such that $\deg g = d_1+d_2+\ldots + d_n$ and such that $g$ has a non-zero coefficient at $X_1^{d_1} \cdots X_n^{d_n}$. Then $g$ does not vanish identically on any additive subgroup $H$ of $F^n$ that spans $F^n$ as a vector space over $F$.
\end{thm}
\begin{prfd}
It suffices to apply the previous theorem to $m=1$, taking all $S_i$ to be the prime field $\F_p$ of $F$, and letting $b_1,b_2,\ldots,b_n$ be a basis for $F^n$ that is contained in $H$. This proves the theorem.
\end{prfd}

\begin{exe}
Prove the following characteristic zero version of the previous theorem. Let $F$ be a field of characteristic zero, let $n\in \Z_{\ge 0}$, and let $g \in F[X_1,\ldots, X_n]$ be non-zero. Then $g$ does not vanish identically on any additive subgroup $H$ of $F^n$ that spans $F^n$ as a vector space over $F$.
\end{exe}

\begin{exe}
Let $k$ be a finite field, let $F \supset k$ be a field extension, let $n\in \Z_{\ge 0}$, and let $g\in F[X_1,\ldots,X_n]$. Suppose that there exist integers $d_1,\ldots,d_n \in \{0,1,\ldots,\#k-1\}$ such that $\deg g = d_1 + d_2 + \ldots + d_n$ and such that $g$ has a non-zero coefficient at $X_1^{d_1} \cdots X_n^{d_n}$. Prove that $g$ does not vanish identically on any sub-$k$-vector space of $F^n$ that spans $F^n$ as an $F$-vector space.
\end{exe}

\subsection{The Jacobson radical}\label{sectiejacob}
\noindent The present section is devoted to two properties of the Jacobson radical that are used in our proof of the unit theorem.

Let $A$ be a ring. We write $A^\ast$ for the group of units of $A$. By an $A$-module, we shall always mean a left $A$-module. An $A$-module $M$ is called \textit{simple} if it has exactly two submodules (namely, $\{0\}$ and $M$). The \textit{Jacobson radical} of $A$ is
\[ \J(A) = \{a\in A : aM = 0 \textnormal{ for each simple } A \textnormal{-module } M\}.\]
It is clearly a two-sided ideal of $A$.

\begin{prop}
Let $A$ be a ring, and $a \in A$. Then:
\[ a \in A^\ast \iff a + \J(A) \in (A/\J(A))^\ast.\]
\end{prop}
\begin{prfd}
The implication $\Rightarrow$ is clear. For $\Leftarrow$, define $U = \{a \in A : a + \J(A) \in (A/\J(A))^\ast\}$, and let $a \in U$. Then $Aa + \J(A) = A$. If $Aa \neq A$, then $A$ has a maximal left ideal $I$ with $Aa \subset I$; for any such $I$, the module $A/I$ is simple, so $\J(A)\cdot (A/I) = 0$, and therefore $\J(A) \subset I$, which leads to the contradiction $A = Aa + \J(A) \subset I$. This proves $Aa = A$, so we can choose $b \in A$ with $ba = 1$. Then $b+ \J(A)$ is the inverse of $a + \J(A)$ in the group $(A/\J(A))^\ast$, so $b \in U$. Therefore every element of $U$ has a left inverse in $U$, which implies that $U$ is a group under multiplication, so $U \subset A^\ast$. This proves the proposition.
\end{prfd}

\begin{exe}\label{exercisejac}
Let $A$ be a ring.
\begin{itemize}[noitemsep, nolistsep]
\item[(a)] Let $I$ be a two-sided ideal of $A$. Prove: $I \subset \J(A) \iff 1 + I \subset A^\ast$.

\item[(b)] The \textit{ring} $A^{\opp}$ \textit{opposite to} $A$ has the same additive group as $A$, but the product of two elements $x$ and $y$ in $A^{\opp}$ is defined to be the product of $y$ and $x$ in $A$. Prove: $\J(A^{\opp}) = \J(A)$.
\end{itemize}
\end{exe}

\noindent We call a ring $D$ a \textit{division ring} if $D^\ast = D \setminus \{0\}$. By a \textit{division algebra} over a commutative ring we mean an algebra that is a division ring. For a ring $D$ and $m\in \Z_{\ge 0}$, we write $\M(m,D)$ for the ring of $m\times m$-matrices with entries in $D$.

\begin{thm}
Let $F$ be a field, and let $A$ be a finite-dimensional algebra over $F$. Then for some $t\in \Z_{\ge 0}$ there exist positive integers $m_1$, $\ldots$, $m_t$ and finite-dimensional division algebras $D_1$, $\ldots$, $D_t$ over $F$ such that 
\[ A/\J(A) \cong \prod_{i=1}^t \M(m_i,D_i) \] as $F$-algebras.
\end{thm}
\begin{prfd}
We start by proving that, as an $A$-module, $A/\J(A)$ is the direct sum of a finite collection of simple $A$-modules. To this end, define $\Is$ to be the set of those left ideals $I$ of $A$ for which $A/I$ can be written as such a direct sum; for example, one has $A \in \Is$. Since $A$ and all of its left ideals are finite-dimensional $F$-vector spaces, we can choose $I \in \Is$ minimal under inclusion. We claim $I = \J(A)$. Write $A/I = \bigoplus_{M \in \Ss} M$, where $\Ss$ is a finite set of simple $A$-modules. Since $\J(A)\cdot M = 0$ for each $M \in \Ss$, we have $\J(A)\cdot (A/I) = 0$, so $\J(A) \subset I$. Assume, for a contradiction, that the opposite inclusion is not valid. Then we can choose a simple $A$-module $N$ and $x\in N$ such that $Ix \neq 0$. By simplicity of $N$, we have $Ix = N$, so the $A$-linear map $f\colon A \lra (A/I)\oplus N, f(a) = (a+I,ax)$, fits in a commutative diagram
\[\xymatrix{ 0 \ar[r] & I \ar[r] \ar[d] & A \ar[r] \ar[d]^f & A/I \ar[r] \ar[d]^{\id} & 0 \\ 0 \ar[r] & N \ar[r] & (A/I)\oplus N \ar[r] & A/I \ar[r] & 0}\]
of $A$-modules that has exact rows and in which the first vertical arrow is surjective. Then $f$ is surjective as well, from which it follows that $\ker f \in \Is$; but $\ker f$ is properly contained in $I$, so this contradicts the minimality of $I$. This proves our claim $I = \J(A)$, so $\J(A) \in \Is$, and therefore $A/\J(A) = \bigoplus_{M\in \Ss} M$ for some finite set $\Ss$ of simple $A$-modules.

We now describe the endomorphism ring $\End_A(A/\J(A))$ of the $A$-module $A/\J(A)$ in two ways. On the one hand, since $A/\J(A)$ is a ring, the $A$-linear maps $A/\J(A) \lra A/\J(A)$ coincide with the right multiplications $x \mapsto xa$ by ring elements $a\in A/\J(A)$. Thus, in the notation of Exercise \ref{exercisejac}(b), we obtain an isomorphism $\End_A(A/\J(A)) \cong (A/\J(A))^{\opp}$ of $F$-algebras.

On the other hand, the additive group of $\End_A(A/\J(A))$ may be identified with
\[ \Hom_A(\bigoplus_{N \in \Ss} N, \bigoplus_{M\in \Ss} M) \cong \bigoplus_{M,N\in \Ss} \Hom_A(N,M),\]
and one readily checks that the product in $\End_A(A/\J(A))$ is, in terms of $\bigoplus_{M,N} \Hom_A(N,M)$, simply given by matrix multiplication. If $f\colon N \lra M$ is a non-zero $A$-linear map, with $N,M \in \Ss$, then by $\ker f \neq N$ and $fN \neq 0$ and simplicity we have $\ker f = 0$ and $fN = M$, so $f$ is an \textit{isomorphism}. Hence $\Hom_A(N,M) = 0$ if $N\not\cong M$, and $\Hom_A(M,M) = \End_A(M)$ is a division algebra; it is finite-dimensional over $F$, since $M$ is. Now if we sort the $M \in \Ss$ by isomorphism, and let $t$ be the number of isomorphism classes occurring, we have $\bigoplus_{M \in \Ss} M \cong \bigoplus_{i = 1}^t M_i^{m_i}$ for certain $m_i \in \Z_{> 0}$, with endomorphism ring isomorphic to $\prod_{i=1}^t \M(m_i,\End_A(M_i))$ as an $F$-algebra.

Comparing the two descriptions of $\End_A(A/\J(A))$, we obtain an isomorphism
\[ A/\J(A) \cong \mleft( \prod_{i=1}^t \M(m_i,\End_A(M_i))\mright)^{\opp} \cong \prod_{i=1}^t \M(m_i,D_i)\]
of $F$-algebras, where $D_i = \End_A(M_i)^{\opp}$ and where the second isomorphism is given by transposition of matrices.

This proves the theorem.
\end{prfd}

\noindent \textbf{Remark.} If $F$ is algebraically closed, then the $D_i$ occurring in the theorem are necessarily equal to $F$, because every finite-dimensional division algebra $D$ over $F$ can be written as a union $\bigcup_{a \in D} F(a)$; each $F(a)$ is commutative, so is a finite field extension of $F$, and therefore equal to $F$ itself.

\subsection{Proof of the unit theorem}
\begin{lem}
Let $F$ be an algebraically closed field, and let $A$ be a finite-dimensional algebra over $F$. Then there exist, for some index set $J$, a basis $(b_i)_{i\in J}$ for $A$ as an $F$-vector space and a polynomial $g\in F[X_i : i \in J]$ such that:
\begin{itemize}
\item[\rm(a)] if $x = \sum_{i \in J} x_i b_i \in A$, with each $x_i \in F$, then $x$ is a unit of $A$ if and only if $g((x_i)_{i \in J}) \neq 0$;

\item[\rm(b)] there exists $m \in \Z_{\ge 0}$ such that each non-zero term of $g$ is of the form $\pm \prod_{i\in I} X_i$ for some subset $I\subset J$ with $\# I = m$.
\end{itemize}
\end{lem}
\begin{prfd}
First suppose that $A$ equals the $m\times m$-matrix algebra over $F$, for some $m \in \Z_{> 0}$. For $i,j \in \{1,\ldots,m\}$, let $b_{ij}\in A$ be the matrix with $(i,j)$-entry $1$ and all other entries $0$. Then $(b_{ij})_{i,j = 1}^m$ is a basis for $A$ over $F$. Also, for $x_{ij} \in F$ and $x = \sum_{i,j} x_{ij} b_{ij}$, we have $x \in A^\ast \iff \det((x_{ij})_{i,j=1}^m) \neq 0$, so we can take $g= \det((X_{ij})_{i,j=1}^m)$; it is well known that each non-zero term of $g$ has the form $\pm \prod_{i=1}^m X_{i\sigma(i)}$ for some permutation $\sigma$ of $\{1,2,\ldots,m\}$. 

Next suppose that $A = A' \times A''$ for certain $F$-algebras $A'$ and $A''$ for which the assertion of the lemma is valid, with bases $(b_i')_{i\in J'}$ and $(b_i'')_{i\in J''}$ and polynomials $g',g''$. We may assume that $J'$ and $J''$ are disjoint. Since $A^\ast = A'^\ast \times A''^\ast$, we may now take the basis for $A$ to be the union of $((b_i',0))_{i\in J'}$ and $((0,b_i''))_{i\in J''}$, with polynomial $g=g'\cdot g''$. 

We pass to the general case. Write $\J(A)$ for the Jacobson radical of $A$. By the theorem of \S 3 and the remark following that theorem, the $F$-algebra $A/\J(A)$ is isomorphic to a product of finitely many matrix algebras over $F$. Hence, by the cases already proved, the lemma is valid for $A/\J(A)$, so there are a basis $(b_i+\J(A))_{i\in J'}$ of $A/\J(A)$ over $F$ and a polynomial $g'\in F[X_i : i \in J']$ satisfying the analogues of (a) and (b). Supplementing $(b_i)_{i\in J'}$ with an $F$-basis for $\J(A)$ we obtain a basis $(b_i)_{i\in J}$ for $A$ over $F$, with $J' \subset J$. By the proposition of \S 3, we can now simply take $g= g'$, viewed as an element of $F[X_i : i \in J]$. This proves the lemma.
\end{prfd}

\begin{exe}
Let the notation and the hypotheses be as in the lemma. Also, let $\Ss$ be a set of simple $A$-modules such that every simple $A$-module is isomorphic to exactly one member of $\Ss$. 
\begin{itemize}[noitemsep,nolistsep]
\item[(a)] Prove that, when $(b_i)_{i\in J}$ and $g$ are constructed as in the proof given, the number $m$ in (b) is given by $m=\sum_{M \in \Ss} \dim_F M$.

\item[(b)] (For readers who know Hilbert's Nullstellensatz.) Prove that, whenever $(b_i)_{i\in J}$ and $g$ satisfy (a) and (b) of the lemma, the number $m$ in (b) is necessarily equal to $\sum_{M \in \Ss} \dim_F M$.
\end{itemize}
\end{exe}

\noindent We can now prove the unit theorem from the introduction, which reads as follows.

\begin{thm}
Let $F$ be a field, let $A$ be a finite-dimensional algebra over $F$, and let $H$ be an additive subgroup of $A$ that spans $A$ as an $F$-vector space. Then $H$ contains a unit of $A$.
\end{thm} 
\begin{prfd}
We proved the theorem already in the introduction if $\Char F = 0$, so we may assume $\Char F = p > 0$.

Let it first be assumed that $F$ is algebraically closed. Choose $(b_i)_{i\in J}$ and $g$ as in the lemma. Since $A^\ast$ is non-empty, part (a) of the lemma implies $g\neq 0$. By part (b), we can choose a non-zero term $\pm \prod_{i \in J} X_i^{d_i}$ of $g$ with $\sum_i d_i = \deg g$ and all $d_i \in \{0,1\}$. The non-vanishing theorem on subgroups, from \S 2, now shows that $g$ does not vanish identically on $H$, so by part (a) of the lemma we have $H\cap A^\ast \neq \emptyset$, as required.

We pass to the general case. Let $\Omega$ be any algebraically closed field extension of $F$, and write $A_{\Omega}$ for the $\Omega$-algebra $A \otimes_F \Omega$. We may identify $A$ with the sub-$F$-algebra $A\otimes 1$ of $A_\Omega$. As an $\Omega$-vector space, $A_\Omega$ is spanned by $A$, and therefore also by $H$. We claim that we have $A^\ast = A \cap A_\Omega^\ast$. Namely, if $a\in A$, then one has $a\in A^\ast$ if and only if the matrix that describes the $F$-linear map $A\lra A$, $x\mapsto xa$, on some $F$-basis for $A$ is invertible; but this matrix doesn't change if we use the same basis for $A_\Omega$ over $\Omega$, so its invertibility is also equivalent to $a \in A_\Omega^\ast$, as required. 

We can now apply the case already proved to $A_\Omega$, and find that $H\cap A_\Omega^\ast$ is non-empty; but $H\subset A$, so we have $H \cap A^\ast = H\cap A \cap A_\Omega^\ast = H\cap A_\Omega^\ast\neq \emptyset$.

This proves the unit theorem.
\end{prfd}

\begin{exe}
Let $F$ be a \textit{perfect} field of characteristic $p>0$, let $A$ be a finite-dimensional algebra over $F$, and put $n = \dim_F A$. Prove that for every basis $(b_i)_{i=1}^n$ of $A$ over $F$ there exists a polynomial $g\in F[X_1,\ldots,X_n]$ that satisfies the hypothesis of the non-vanishing theorem on subgroups from \S 2, such that for any $x = \sum_{i=1}^n x_i b_i \in A$, with $x_i \in F$, one has: $x \in A^\ast \iff g(x_1,\ldots,x_n) \neq 0$.
\end{exe}

\begin{exe}
Let $F$ be a field, let $A$ be a finite-dimensional algebra over $F$, and let $H$ be an additive subgroup of $A$ that spans $A$ as an $F$-vector space. Prove that every coset of $H$ in $A$ contains a unit of $A$.
\end{exe}

\begin{exe}
Let $F$ be a field, let $m\in \Z_{>0}$, and let $H$ be an additive subgroup of the ring $\M(m,F)$ of $m\times m$-matrices over $F$ that spans $\M(m,F)$ as an $F$-vector space. Prove that for every $c\in F$ and every $a \in \M(m,F)$ there exists $x \in a+ H$ with $\det x \neq c$.
\end{exe}

\subsection{Descending module properties}
Throughout this section we make the following assumptions. By $F$ we denote a field, and $K$ is an extension field of $F$. For any $F$-vector space $V$, we shall write $V_K = V \otimes_F K$, which is a $K$-vector space with $\dim_K(V_K) = \dim_F(V)$; we shall always identify $V$ with its image $V\otimes 1$ in $V_K$. We shall constantly use that a sequence $V \lra V' \lra V''$ of $F$-vector spaces is exact if and only if the sequence $V_K \lra V_K'\lra V_K''$ of $K$-vector spaces it induces is exact; this is because $K$ is $F$-free of non-zero rank. In particular, one has $V = 0 \iff V_K = 0$, and an $F$-linear map $V \lra V'$ is zero, injective, surjective, or bijective, if and only if the induced $K$-linear map $V_K \lra V_K'$ has the same respective property.

By $A$ we denote an $F$-algebra. Then $A_K$ is a $K$-algebra and if $M$ is an $A$-module, then $M_K$ is an $A_K$-module. \textbf{Modules denoted by $M$ or $N$ will always be assumed to be finite-dimensional as vector spaces over $F$}; then $M_K$ and $N_K$ have the same property over $K$. We do \textit{not} assume that $A$ be finite-dimensional over $F$; however, when dealing with an $A$-module $M$, we can typically reduce to that case by replacing $A$ by its image in $\End_F(M)$, which does have finite $F$-dimension.

In the present section we are interested in properties of an $A$-module $M$ that are \textit{base-change invariant} in the sense that they pass to the induced $A_K$-module $M_K$. Each such property gives rise to a \textit{descent problem}: if the $A_K$-module $M_K$ has the property, does the same hold for the $A$-module $M$? In many cases where the answer is affirmative, an easy proof can be given by means of the unit theorem; often a slighty sharper statement can be obtained, as was the case with the normal basis theorem in \S 1.

The following lemma will be constantly used. 
\begin{lem}
Let $M$, $N$ be $A$-modules. Then the natural map \[\Hom_A(M,N)_K \lra \Hom_{A_K}(M_K,N_K)\] is an isomorphism.
\end{lem}
\begin{prfd}
Replacing $A$ by its image in $(\End_F M) \times \End_F N$, we may assume $\dim_F(A) < \infty$. For $M=A$, both groups in the lemma may be identified with $N_K$, and the assertion is easy to check. For general $M$, we can now use a finite presentation $A^m \lra A^n \lra M \lra 0$ of $M$ as an $A$-module, and left exactness of $\Hom$, to reduce to the previous case. This proves the lemma.
\end{prfd}

\noindent The first property we consider is \textit{cyclicity}. We call an $A$-module $M$ \textit{cyclic} if there exists $x \in M$ with $M = Ax$ or, equivalently, if there is a surjective $A$-linear map $A\lra M$.

\begin{thm}[Generator theorem]
Let $M$ be a cyclic $A$-module, and let $H$ be an additive subgroup of $M$ that spans $M$ as an $F$-vector space. Then there exists $x \in H$ with $M = Ax$.
\end{thm}
\begin{prfd}
As above, we reduce to the case $\dim_F A < \infty$. Let $f\colon A \lra M$ be a surjective $A$-linear map. Then $f^{-1}H$ is an additive subgroup of $A$ that spans $A$ as an $F$-vector space, so by the unit theorem there is a unit $u \in f^{-1}H$. Then $x = f(u) \in M$ satisfies $Ax = f(Au) = f(A) = M$, as required. This proves the theorem.
\end{prfd}

\begin{rem}
See Exercise 18 for a generalization to modules that require more generators.
\end{rem}

\begin{thm}[Descending cyclicity]
Let $M$ be an $A$-module. Then $M$ is cyclic if and only if $M_K$ is a cyclic $A_K$-module, and one has $A\cong M$ as $A$-modules if and only if $A_K \cong M_K$ as $A_K$-modules.
\end{thm}
\begin{prfd}
Both ``only if'' statements are clear. For the first ``if'' statement, assume that $M_K$ is a cyclic $A_K$-module. Applying the generator theorem, with $M_K$, $A_K$, $K$, $M=M\otimes 1$ in the roles of $M$, $A$, $F$, $H$, we find $x\in M$ such that  the map $A\lra M$, $a\mapsto ax$, becomes surjective after tensoring with $K$. Then it is already surjective, as required. For the second ``if'' statement, the map $A\lra M$ obtained in the same way is also injective, because $A_K \cong M_K$ implies $\dim_F A = \dim_K (A_K) = \dim_K(M_K) = \dim_F M$. This proves the theorem.
\end{prfd}

\noindent Changing rings, we can extend the second part of the theorem above to more general isomorphisms.

\begin{thm}[Descending isomorphisms]
Let $M$, $N$ be $A$-modules. Then one has $M\cong_A N \iff M_K \cong_{A_K} N_K$.
\end{thm}
\begin{prfd}
It suffices to prove $\Leftarrow$. Let $f\colon M_K \lra N_K$ be an $A_K$-linear isomorphism. Then $\End_{A_K}(N_K) \lra \Hom_{A_K}(M_K,N_K)$, $g\mapsto gf$, is an isomorphism of $\End_{A_K}(N_K)$-modules, so by the lemma the $\End_A(N)$-module $\Hom_A(M,N)$ becomes isomorphic to $\End_A(N)$ after tensoring with $K$. Applying the theorem proved above, with $\End_A N$ and $\Hom_A(M,N)$ in the roles of $A$ and $M$, one finds $\End_A N \cong\Hom_A(M,N)$ as $\End_A N$-modules, so $\Hom_A(M,N) = (\End_A N) \cdot h$ for some $A$-linear map $h\colon M \lra N$. Then one also has \linebreak $\Hom_{A_K}(M_K,N_K) = \End_{A_K}(N_K)\cdot h_K$ (for $h_K = h\otimes \id_K$), so $f = i\cdot h_K$ for some $i \in \End_{A_K}(N_K)$. Since $f$ is an isomorphism $h_K$ is injective, and comparing dimensions one finds that $h_K$ is an isomorphism. Then $h$ is an isomorphism, as desired.
\end{prfd}

\begin{exe}
Suppose that $M$ and $N$ are isomorphic $A$-modules. Let $H$ be a subgroup of $\Hom_A(M,N)$ that spans $\Hom_A(M,N)$ as an $F$-vector space. Prove that $H$ contains an isomorphism $M \lra N$ of $A$-modules.
\end{exe}

\begin{thm}[Noether--Deuring]
Let $M$ and $N$ be $A$-modules. Then there exists an $A$-module $P$ with $M \cong N \oplus P$ as $A$-modules if and only if there exists an $A_K$-module $Q$ with $M_K \cong N_K \oplus Q$ as $A_K$-modules.
\end{thm}
\begin{prfd}
For the ``only if''-part one can take $Q = P_K$. For the ``if''-part, we write $V = \Hom_A(M,N)$ and $W = \Hom_A(N,M)$. Note that by the lemma, we may identify $V_K$ with $\Hom_{A_K}(M_K,N_K)$ and $W_K$ with $\Hom_{A_K}(N_K,M_K)$. Likewise, if we write $B = \End_A(N)$, which is a finite-dimensional algebra over $F$, then we may identify $B_K$ with $\End_{A_K} N_K$. It is an easy exercise to show that the existence of $P$ as in the theorem is equivalent to the existence of a pair $(f,g) \in V \times W$ with $f\cdot g \in B^\ast$.  So, to prove the theorem, it suffices to prove that such a pair exists if we are given a pair $(f',g') \in V_K \times W_K$ with $f'\cdot g' \in B_K^\ast$. The $K$-span of the subgroup $V\cdot g' \subset B_K$ equals $V_K\cdot g'$, which is a left $B_K$-ideal containing the unit $f'\cdot g'$ and is therefore equal to $B_K$. Hence, applying the unit theorem with $K$, $B_K$, $V\cdot g'$ in the respective roles of $F$, $A$, $H$, we find $f\in V$ with $f\cdot g' \in B_K^\ast$. Next, the $K$-span of $f\cdot W \subset B_K$ is a right $B_K$-ideal containing the unit $f\cdot g'$, so it equals $B_K$, and with $K$, $B_K$, $f\cdot W$ in the roles of $F$, $A$, $H$ we find $g\in W$ with $f\cdot g \in B_K^\ast$. Since also $f\cdot g \in B$, we have $f\cdot g \in B^\ast$, as desired. This proves the theorem.
\end{prfd}

\noindent In the exercises that follow, the assumptions stated at the beginning of this section remain in force. This applies in particular to the assumption $\dim_F M < \infty$. We recall that a module $P$ over a ring $R$ is called \textit{projective} if for any surjective $R$-linear map $V\lra W$ of $R$-modules the induced map $\Hom_R(P,V) \lra \Hom_R(P,W)$ of abelian groups is surjective; equivalently, $P$ is projective if and only if there exists an $R$-module $Q$ such that $P \oplus Q$ is a \textit{free} $R$-module.
\vspace{0.5\baselineskip}

\begin{exe}
Let $R$ be a ring, let $P$ be a projective $R$-module, and let $V$ be an $R$-module. Prove that the following two assertions are equivalent:
\begin{itemize}[noitemsep,nolistsep]
\item[(i)] there is a surjective $R$-linear map $f \colon P \lra V$;

\item[(ii)] there exists $f \in \Hom_R(P,V)$ with $f\cdot \End_R P = \Hom_R (P,V)$, and every such $f$ is surjective.
\end{itemize}
\end{exe}

\begin{exe}
Let $n\in \Z_{\ge 0}$, and let $M$ be an $A$-module that can be generated by $n$ elements. Prove that for every additive subgroup $H\subset M^n$ that spans $M^n$ as an $F$-vector space there exists $(x_1,x_2,\ldots,x_n) \in H$ such that $M = \sum_{i=1}^n Ax_i$.
\end{exe}

\begin{exe}
Let $M$ be an $A$-module, and let $P$ be a finitely generated projective $A$-module. Prove: there exists a surjective $A$-linear map $P\lra M$ if and only if there exists a surjective $A_K$-linear map $P_K\lra M_K$.
\end{exe}

\begin{exe}
Let $M$ be an $A$-module. Prove that the least number of generators of $M$ as an $A$-module is equal to the least number of generators of $M_K$ as an $A_K$-module.
\end{exe}

\begin{exe}
Suppose that $F$ is a \textit{finite} field, let $V$ be a $2$-dimensional $F$-vector space, and write $\P$ for the set of $1$-dimensional subspaces of $V$.
\begin{itemize}[noitemsep,nolistsep]
\item[(a)] Put $A = F\oplus F^\P$, and define a multiplication on $A$ by $(a,x) \cdot (b,y) = (ab, ay + bx)$. Prove that $A$ is a commutative ring and that the map $F\lra A$, $a\mapsto (a,0)$ is a ring homomorphism that makes $A$ into a finite-dimensional algebra over $F$.

\item[(b)] Put $M = V \oplus \bigoplus_{W\in \P}(V/W)$. Prove that the scalar multiplication $\cdot\colon A\times M \lra M$, $(a,(b_W)_{W\in \P})\cdot (x,(y_W + W)_{W\in \P}) = (ax, (ay_W + b_W x + W)_{W\in \P})$ makes $M$ into an $A$-module satisfying $\dim_F M < \infty$.

\item[(c)] Prove: there exists an injective $A_K$-linear map $A_K\lra M_K$ if and only if $K \neq F$.

\end{itemize}
\end{exe}

\begin{rem}
Since $A$ is projective as an $A$-module, Exercise 21(c) shows that Exercise 19 may fail when injective instead of surjective maps are considered. The following exercise shows that Exercise 19 also becomes incorrect when the projectivity condition is omitted.
\end{rem}

\begin{exe}
Let $F$, $A$ be as in Exercise 21. Construct $A$-modules $P$ and $M$, both finite-dimensional over $F$, such that there exists a surjective $A_K$-linear map $P_K \lra M_K$ if and only if $K\neq F$.
\end{exe}

\subsection{Glynn's determinant formula}
In the present section we give a new proof of the commutative case of the unit theorem. Throughout this section, $p$ is a prime number, $F$ is a field of characteristic $p$, and $n\in \Z_{\ge 0}$.

\begin{thm}[Glynn's determinant formula, 1998]
For any $n\times n$-matrix $A = (a_{ij})_{i,j=1}^n$ over $F$, the coefficient of the polynomial 
\[\prod_{j=1}^n\mleft(\sum_{i=1}^n a_{ij} X_i\mright)^{p-1} \textit{ at } X_1^{p-1} X_2^{p-1}\cdots X_n^{p-1}\]
equals $(\det A)^{p-1}$.
\end{thm}

\begin{exam}
If $A$ is a \textit{diagonal} matrix, then an easy computation shows that one has in fact \[ \prod_{j=1}^n \mleft( \sum_{i=1}^n a_{ij}X_i\mright)^{p-1} = (\det A)^{p-1}\cdot X_1^{p-1}\cdots X_n^{p-1}.\]
\end{exam}

\begin{prfd}[Proof of the theorem]
Write $V$ for the $F$-vector space  $\sum_{i=1}^n F\cdot X_i$. Then $F[X_1,X_2,\ldots,X_n]$ may be identified with the symmetric algebra of $V$ over $F$. Hence the group $\Aut V$, which may be identified with $\GL(n,F)$, acts on the $F$-algebra $F[X_1,X_2,\ldots,X_n]$ in a natural way. Write $I$ for the ideal generated by $\{v^p:v\in V\}$. Then for all $\sigma \in \Aut V$ one has $\sigma I = I$. Also, since $p$-th powering is additive, the ideal $I$ is generated by $\{X_1^p,\ldots, X_n^p\}$ and has therefore as an $F$-basis the set of monomials $X_1^{d_1}\cdots X_n^{d_n}$ in which at least one $d_i$ is at least $p$. It follows that the $F$-algebra $R = F[X_1,\ldots,X_n]/I$ has the images of the monomials $X_1^{d_1}\cdots X_n^{d_n}$ with all $d_i < p$ as an $F$-basis. Hence $R$ contains $V$, and the $R$-ideal $J = \{r\in R : rV = 0\}$ is the one-dimensional $F$-vector space generated by the image of $X_1^{p-1}\cdots X_n^{p-1}$ in $R$.

The natural action of $\Aut V$ on $R$ preserves $J$ so yields a group homomorphism $\varphi\colon \Aut V \lra \Aut_F J \cong F^\ast$. The latter group is abelian, so $\ker \varphi$ contains the commutator subgroup $[\Aut V, \Aut V]$. But if $\# F > 2$---which we may assume, upon replacing $F$ by a suitable extension---then by Exercise 23 we have $[\Aut V, \Aut V] = \ker \det$, where $\det\colon \Aut V \lra F^\ast$ is the determinant map. Therefore, for $\sigma \in \Aut V$, the image $\varphi(\sigma)$ depends only on $\det \sigma$, and checking the case that $\sigma$ is given by a diagonal matrix, one finds $\varphi(\sigma) = (\det \sigma)^{p-1}$. Rewriting this in terms of the $X_i$, one obtains from this Glynn's formula in the case $A$ is invertible. In the generic case, when the $a_{ij}$ are $n^2$ polynomial variables over $\F_p$ and $F = \F_p(a_{ij} : 1 \le i,j \le n)$, the matrix is indeed invertible; and the general case follows from the generic case. This proves the theorem.
\end{prfd}

\begin{exe}
Prove: if $V$ is an $n$-dimensional vector space over $F$, then $[\Aut V, \Aut V]$ equals the kernel of $\det\colon \Aut V \lra F^\ast$, except in the case $n = \# F = 2$, in which case $[\Aut V, \Aut V]$ has index $2$ in $\ker \det$.
\end{exe}

\begin{exe}
Let $m\in \Z_{\ge 0}$. Prove that Glynn's determinant formula is also valid when one replaces $p-1$, at all occurrences, by $p^m-1$.
\end{exe}

\begin{exe}
Let $k$ be a finite field, and $d_1$, $d_2$, $\ldots$, $d_n \in \Z_{\ge 0}$. Prove:
\[ \sum_{x_1, x_2, \ldots, x_n \in k} x_1^{d_1}\cdots x_n^{d_n} = \begin{cases} (-1)^n & \textnormal{ if each $d_i$ is a positive multiple of $\# k -1$,} \\ 0 & \textnormal{ otherwise.}\end{cases}\]
\end{exe}
\vspace{\baselineskip}
\begin{thm}
Let $B = \{(b_{ij})_{j=1}^n : 1 \le i \le n\}$ be a basis of the $F$-vector space $F^n$, and let $H\subset F^n$ be the additive subgroup generated by $B$. Then we have 
\[ \sum_{(y_1,y_2,\ldots,y_n) \in H} y_1^{p-1}y_2^{p-1}\cdots y_n^{p-1} = (-1)^n \cdot (\det(b_{ij})_{i,j=1}^n)^{p-1}.\]
\end{thm}
\begin{prfd}
If we write $f = \prod_{j=1}^n \mleft(\sum_{i=1}^n b_{ij} X_i\mright)^{p-1}$, then the sum on the left equals
\[ \sum_{x_1,x_2,\ldots,x_n \in \F_p} f(x_1,\ldots,x_n)^{p-1}.\] We apply Exercise 25 to $k = \F_p$. Since $f$ is homogeneous of degree $n\cdot (p-1)$, it has at most one non-zero monomial $c\cdot X_1^{d_1}\cdots X_n^{d_n}$ with all $d_i \in \Z_{> 0}\cdot (p-1)$, namely with all $d_i = p-1$; and the coefficient $c$ at that monomial equals $(\det(b_{ij})_{i,j=1}^n)^{p-1}$, by Glynn's formula. Hence the theorem follows from Exercise 25. This proves the theorem.
\end{prfd}

\noindent Since the right hand side in the formula just proved is non-zero, it implies that at least one $(y_1,y_2,\ldots,y_n)\in H$ satisfies $y_1y_2\cdots y_n\neq 0$. As we saw in section $1$, this implies the commutative case of the unit theorem.

\begin{exe}
Let $B = \{ (b_{ij})_{j=1}^n : 1 \le i \le n\}$ be a basis of the $F$-vector space $F^n$, and let $k$ be a \textit{finite} subfield of $F$. Write $W$ for the $k$-vector space generated by $B$. Prove:
\[ \sum_{(y_1,y_2,\ldots, y_n)\in W} y_1^{\# k-1} y_2^{\# k - 1} \cdots y_n^{\#k -1} = (-1)^n\cdot (\det(b_{ij})_{i,j=1}^n)^{\#k-1}.\]
\end{exe}

\subsection{A lower bound and an algorithm}
We give a third proof of the unit theorem in the commutative case, again concentrating on totally split algebras. It will lead to a lower bound, as stated in the theorem at the end of the introduction. 

Throughout this section, $F$ is a field.

\begin{thm}
Let $E$ be a subfield of $F$, and $n\in \Z_{\ge 0}$. Let $W \subset F^n$ be a sub-$E$-vector space of $F^n$ such that the natural map $W\otimes_E F \lra F^n$ is an isomorphism. Then $W$ has a basis $B$ over $E$ such that $\sum_{b\in B} E^\ast \cdot b \subset F^{\ast n}$.
\end{thm}
\begin{prf}
In this proof, we abbreviate $-\otimes_E F$ to $-_F$. For $i \in \{1,2,\ldots, n\}$, let $\pi_i\colon W_F \lra F$ be the composition of the isomorphism $W_F \stackrel{\sim}{\lra} F^n$ and the $i$-th coordinate projection $F^n \lra F$. Then $W_F \lra F^n$, $x \mapsto (\pi_i(x))_{i=1}^n$, is an isomorphism of $F$-vector spaces, which is the same as saying that $\pi_1$, $\pi_2$, $\ldots$, $\pi_n$ form a basis for the $F$-vector space $\Hom_F(W_F,F)$ or, equivalently, that $\pi_1\wedge \pi_2 \wedge \ldots \wedge \pi_n$ is a non-zero element of the exterior power $\bigwedge^n \Hom_F(W_F,F)$, which is a $1$-dimensional $F$-vector space.

For each $i$, define $W_i = W \cap (\ker \pi_i)$. The natural map $\Hom_E(W/W_i,E)_{F} \lra \Hom_F(W_F/(W_i)_F,F)$ is an isomorphism (cf.\ the lemma in section 5). Here we may identify $\Hom_F(W_F/(W_i)_F,F)$ with the set of all $f\in \Hom_F(W_F,F)$ that are $0$ on $(W_i)_F$; so $\pi_i \in \Hom_F(W_F/(W_i)_F,F)$. The image of the natural map
\[ \mleft(\bigotimes_{i=1}^n \Hom_E(W/W_i,E)\mright)_{\!\!\!F} \stackrel{\sim}{\to} \bigotimes_{i=1}^n \Hom_F(W_F/(W_i)_F,F) \to \bigwedge^n \Hom_F(W_F,F)\]
contains the non-zero element $\pi_1 \wedge \pi_2 \wedge \ldots \wedge \pi_n$. It follows that there are $\rho_i \in \Hom_E(W/W_i,E)$ such that the image of $\rho_1\otimes \rho_2\otimes \ldots \otimes \rho_n$ is non-zero. Then $\rho_1\wedge \rho_2 \wedge \ldots \wedge \rho_n$ is non-zero as an element of $\bigwedge^n \Hom_E(W,E)$, so the $E$-linear map $W \lra E^n$, $x\mapsto (\rho_i(x))_{i=1}^n$, is an isomorphism. We now choose $B$ to be the inverse image of the standard basis of $E^n$. Then for each $x \in \sum_{b\in B} E^\ast \cdot b$ and each $i$ one has $\rho_i(x) \neq 0$, so $x\notin W_i = W \cap \ker \pi_i$, so $\pi_i(x) \neq 0$; and therefore $x \in F^{\ast n}$. This proves the theorem.
\end{prf}

\begin{thm}[Lower bound]
Let $F$ be a field, let $A$ be a commutative algebra over $F$ of finite vector space dimension $n$, and let $H\subset A$ be an additive subgroup that spans $A$ as an $F$-vector space. Assume that the characteristic $p$ of $F$ is positive. Then $\# (H\cap A^\ast) \ge (p-1)^n$.
\end{thm}
\begin{prfd}
As we argued in the introduction, we may assume that $A$ is totally split, so $A \cong F^n$ as $F$-algebras, and that $H$ is the $\F_p$-span of some $F$-basis for $A$. With $E = \F_p$ and $W= H$ we can now apply the previous theorem, and we find an $\F_p$-basis $B$ for $H$ such that $\sum_{b\in B} \F_p^\ast \cdot b \subset A^\ast$.~This~implies~the~theorem.
\end{prfd}

\noindent For every $F$ of non-zero characteristic $p$ and every $n\in \Z_{\ge 0}$, one readily constructs $A$ and $H$ such that equality holds in the theorem. It is of interest to classify all cases in which the equality sign is valid, and also to drop the commutativity assumption.

\begin{alg}
The proofs of the unit theorem that we gave in the earlier sections, in the general and in the commutative case, do not offer an efficient method for finding the unit that is asserted to exist. By contrast, the proof that we gave in the present section does lead to an efficient algorithm in the totally split case. It is efficient in the sense that, when sufficiently good algorithms are available for doing arithmetic in the fields involved, it runs in polynomial time.

Resuming the notation used in the first theorem above and its proof, we outline an algorithm that, given a field extension $E \subset F$, a non-negative integer $n$, and a basis $C$ for $F^n$ over $F$, constructs a basis $B$ for $W =\sum_{c\in C} E\cdot c$ over $E$ such that $\sum_{b\in B} E^\ast \cdot b \subset (F^{\ast})^n$. 

The algorithm starts by computing $E$-bases for the subspaces $W_i \subset W$ defined in the proof, as well as bases $\Ph_i$ for the $E$-vector spaces $\Hom_E(W/W_i,E)$, expressed on a basis for $\Hom_E(W,E)$. For $\rho \in \Ph_i$, write $\rho_F$ for $\rho\otimes \id_F$, which we view as an element of $\Hom_F(W_F,F)$. Let now $\pi_i \in \Hom_F(W_F,F)$ be as in the proof, so that $\pi_1\wedge \ldots \wedge \pi_n \neq 0$. The algorithm finds, for $i = 1,2,\ldots, n$ in succession, an element $\rho_i \in \Ph_i$ such that $\rho_{1F}\wedge \rho_{2F} \wedge \ldots \wedge \rho_{iF} \wedge \pi_{i+1} \wedge \ldots \wedge \pi_n \neq 0$.  This is done as follows. Assume, inductively, that $\rho_{1F}\wedge \rho_{2F} \wedge \ldots \wedge \rho_{i-1F} \wedge \pi_{i} \wedge \pi_{i+1} \wedge \ldots \wedge \pi_n \neq 0$. As we saw in the proof, $\pi_i$ belongs to $\Hom_F(W_F/(W_i)_F,F)$, which is the $F$-linear span of $(\Ph_i)_F$, so for at least one $\rho \in \Ph_i$ one has $\rho_{1F} \wedge \ldots \wedge \rho_{i-1F}\wedge \rho_F \wedge \pi_{i+1}\wedge \ldots \wedge \pi_n \neq 0$, or, equivalently, the $F$-linear map $W_F \lra F^n$, \[x\mapsto (\rho_{1F}(x),\ldots, \rho_{i-1F}(x), \rho_F(x), \pi_{i+1}(x),\ldots,\pi_n(x))\] is an isomorphism. The first $\rho \in \Ph_i$ that is found to satisfy the latter condition is chosen to be $\rho_i$. Once $\rho_1$, $\rho_2$, $\ldots$, $\rho_n$ have all been found, one has an isomorphism $W \lra E^n$, $x\mapsto (\rho_i(x))^{n}_{i=1}$ of $E$-vector spaces. As in the proof, one now takes $B$ to be the inverse image of the standard basis of $E^n$.

It would be of great interest to extend the algorithm to the non-commutative case.
\end{alg}

\subsection{References}
\begin{ack} are due to Bas Edixhoven for much of section $7$; to Jan Draisma for telling me about Glynn's formula; to Martin Heemskerk for writing [8] and for drawing my attention to the Noether--Deuring theorem; to Lenny Taelman for various useful remarks; and to Abtien Javanpeykar for \LaTeX\ help. The Max Planck Institute for Mathematics in Bonn is thanked for its hospitality and its support.
\end{ack}

\begin{bib}
Much of the material in sections $1$, $4$ and $5$ is taken from Martin Heemskerk's bachelor thesis [8]. The combinatorial Nullstellensatz is due to Noga Alon [1]. A particularly short proof was given by Terence Tao [11]. For the theorem of Cauchy--Davenport, see [5, 6, 12]. More about the Jacobson radical can be found in [10, 9]. For descent properties as in section $5$, see [3, Ch.\! VIII, \S 2.5, Th.\! 3; 4, Ch.\! V, Annexe, Prop.\! 1; 9, Theorem 19.25]. Glynn's determinant formula appeared in [7]. A solution to Exercise 23 can be found in [2].
\end{bib}

\singlespacing \small
\begin{itemize}
\item[1.] N. Alon, Combinatorial Nullstellensatz, Combin.\! Probab.\! Comput.\! \textbf{8} (1999), 7--29.
\item[2.] E. Artin, \textit{Geometric algebra}, Interscience Publishers, Inc., New York, 1957.
\item[3.] N. Bourbaki, \textit{\'{E}l\'{e}ments de math\'{e}matique, Alg\`{e}bre, Chapitre 8}, Springer, Heidelberg, 2012.
\item[4.] N. Bourbaki, \textit{\'{E}l\'{e}ments de math\'{e}matique, Groupes et alg\`{e}bres de Lie, Chapitres 4, 5 et 6}, Hermann, Paris, 1968.
\item[5.] A. L. Cauchy, \textit{Recherches sur les nombres}, J. \'{E}cole Polytech.\! \textbf{9} (1813), 99--116.
\item[6.] H. Davenport, \textit{On the addition of residue classes}, J.\! London Math.\! Soc.\! \textbf{10} (1935), 30--32.
\item[7.] D. G. Glynn, \textit{The modular counterparts of Cayley's hyperdeterminants}, Bull.\! Austral.\! Math.\! Soc.\! \textbf{57} (1998), 479--492.
\item[8.] M. Heemskerk, \textit{Basisuitbreidingen en de Combinatorische Nullstellensatz}, Ma\-thematisch Instituut, Universiteit Leiden, 2014, see:\\ \textsc{\footnotesize https://www.math.leidenuniv.nl/scripties/HeemskerkBach.pdf}. 
\item[9.] T. Y. Lam, \textit{A first course in noncommutative rings}, Springer--Verlag, New York, 1991.
\item[10.] S. Lang, \textit{Algebra}, revised third edition, Springer--Verlag, New York, 2002.
\item[11.] T. Tao, \textit{Algebraic combinatorial geometry: the polynomial method in arithmetic combinatorics, incidence combinatorics, and number theory}, EMS Surv.\! Math.\! Sci.\! \textbf{1} (2014), 1--46.
\item[12.] T. Tao, V. Vu, \textit{Additive combinatorics}, Cambridge University Press, Cambridge, 2006.

\end{itemize}

\noindent \footnotesize \textsc{Mathematisch Instituut, Universiteit Leiden,\\
Postbus 9512, 2300 RA Leiden, Nederland.}\\
\tt hwl@math.leidenuniv.nl
\end{document}